\documentclass[11pt,twoside]{amsart}
\def\RR{{\mathbb R}}
\def\CC{{\mathbb C}}
\newtheorem{theorem}{Theorem}[section]

\newtheorem{proposition}[theorem]{Proposition}
\numberwithin{equation}{section}

\begin{document}

\title[Homogeneous Tube Domains]{Towards a Classification of Homogeneous Tube
Domains in $\CC^4$}

\author[M. Eastwood]{Michael Eastwood}
\address{Department of Mathematics\\
University of Adelaide\\
SA 5005,\newline
AUSTRALIA}
\email{meastwoo@maths.adelaide.edu.au}
\thanks{The first author is supported by the Australian Research Council}

\author[V. Ezhov]{Vladimir Ezhov}
\address{School of Mathematics and Statistics\\
University of South Australia\\
Mawson Lakes Blvd\\
Mawson Lakes\\
SA 5091\\
AUSTRALIA}
\email{vladimir.ejov@unisa.edu.au}
\thanks{The second author gratefully acknowledges ARC Discovery grant
DP0450725}

\author[A. Isaev]{Alexander Isaev}
\address{Department of Mathematics\\
The Australian National University,\newline
Canberra\\
ACT 0200\\
AUSTRALIA}
\email{alexander.isaev@maths.anu.edu.au}

\renewcommand{\subjclassname}{\textup{2000} Mathematics Subject Classification}
\subjclass{Primary 32M10; Secondary 32V15}
\renewcommand{\subjclassname}{\textup{2000} Mathematics Subject Classification}
\keywords{Unbounded domains, Holomorphic homogeneity}

\begin{abstract}
We classify the tube domains in $\CC^4$ with affinely homogeneous base whose
boundary contains a non-degenerate affinely homogeneous hypersurface. It
follows that these domains are holomorphically homogeneous and amongst them
there are four new examples of unbounded homogeneous domains (that do not have
bounded realisations). These domains lie to either side of a pair of
Levi-indefinite hypersurface. Using the geometry of these two
hypersurfaces, we find the
automorphism groups of the domains.
\end{abstract}
\maketitle

\setcounter{section}{-1}
\section{Introduction}

The study of holomorphically homogeneous domains in complex space goes back to
\'E. Cartan \cite{C} who determined all bounded symmetric domains in $\CC^n$ as
well as all bounded homogeneous domains in $\CC^2$ and $\CC^3$. Due to the
fundamental theorem of Vinberg, Gindikin, and Pyatetskii-Shapiro, every bounded
homogeneous domains can be realised as a Siegel domain of the second kind (see
\cite{P-S}). Although this result does not immediately imply a complete
classification of bounded homogeneous domains, it reduces the classification
problem to that for domains of a very special form. A generalisation of the
above theorem to the case of unbounded domains for the class of rational
homogeneous domains was obtained by Penney in his remarkable paper \cite{P2},
where the role of models is played by so-called Siegel domains of type $N$-$P$.
Nevertheless, the classification problem for the unbounded case is far from
fully understood.

In \cite{L2}, Loboda shows that any holomorphically
homogeneous non-spherical tube hypersurface in $\CC^2$ has an
affinely homogeneous
base. The generalisation of this result to higher dimensions is an open
problem. \mbox{Loboda}'s result motivates our study. We shall consider tubes
over affinely homogeneous bases and suppose that the boundary
of the base contains an affinely homogeneous hypersurface. In these
circumstances, we can classify the resulting tube domains in low dimensions and
determine their automorphism groups. The first interesting examples occur
in~$\CC^4$.

A simple natural class of unbounded homogeneous domains comes from generalising
the unbounded form of the unit ball
$$B^n:=\left\{z\in\CC^n: \hbox{Re}\,z_n>|z'|^2\right\}$$
and the complementary domain $\left\{z\in\CC^n: \hbox{Re}\,z_n<|z'|^2\right\}$,
where $z=(z_1,\dots,z_n)$ are coordinates in $\CC^n$ and
$z':=(z_1,\dots,z_{n-1})$. Indeed, the domains lying on either side of the
quadric given by the equation
\begin{equation}\label{quadric}\hbox{Re}\,z_n=\langle z',z'\rangle,
\end{equation}
where $\langle z',z'\rangle$ is a Hermitian form in $\CC^{n-1}$, are easily
seen to be homogeneous. Among them, $B^n$ is the only domain, up to a linear
change of coordinates, that admits a bounded realisation. Domains of this kind
with $\langle z',z'\rangle$ non-degenerate were generalised by Penney, who
introduced and studied a special class of Siegel $N$-$P$ domains called
(homogeneous) {\em nil-balls}\, (see \cite{P1,P3}). A nil-ball is a certain
$N$-$P$ domain on which a solvable algebraic group acts transitively and
polynomially by biholomorphic transformations, with smooth connected algebraic
Levi non-degenerate boundary on which a codimension 1 nil-radical of the group
acts simply transitively.

The examples arising from Hermitian forms, as above, can also be written as
tube domains, i.e.\ domains of the form $D_{\Omega}:=\Omega+i\RR^n$, where the
base $\Omega$ is a domain in $\RR^n\subset\CC^n$ (see, e.g., \cite{I} for an
explicit change of coordinates). Tube domains are clearly unbounded and often
do not have bounded realisations. If $\Omega$ is affinely homogeneous, then
$D_{\Omega}$ is holomorphically homogeneous since every affine mapping of
$\RR^n$ can be lifted to an affine mapping of $\CC^n$ and since $D_{\Omega}$ is
invariant under translations in imaginary directions. In \cite{P1} a large
class of nil-balls of this kind is introduced.

The new examples of tube domains in $\CC^4$ that we give in this paper have
affinely homogeneous bases,  possess no bounded realisations, and are not
equivalent to nil-balls. They arise by the following
construction/classification. Start with an affinely homogeneous hypersurface
$\Gamma\subset\RR^n$. We might find that the orbits of the affine symmetry
group of $\Gamma$ are, in addition to $\Gamma$ itself, some domain or domains
$\Omega$ to either side of it. Of course, in this case, the action of the
symmetry group on $\Gamma$ must have isotropy. A classification of suitable
$\Gamma$ and $\Omega$ leads to a classification of corresponding homogeneous
tubes~$D_\Omega$.

Affinely homogeneous hypersurfaces $\Gamma\subset\RR^n$ have been classified
for $n=2$ (see, e.g., \cite{NS1}) and $n=3$ in \cite{DKR,EE1}. We are
interested in hypersurfaces with non-trivial affine isotropy subgroup, and when
$n=4$ such $\Gamma$ have also been classified (see \cite{EE2}) under
the additional hypothesis of non-degeneracy (equivalently, under the
hypothesis that $\Gamma+i\RR^n$ is Levi non-degenerate). In this
paper we concentrate on the case of non-degenerate $\Gamma$ for
$n=4$. This leads to a
relatively small number of affinely homogeneous domains: most of the
automorphism groups from~\cite{EE2} do not have open orbits. This, in turn,
leads to the following classification in~$\CC^4$:--
\begin{theorem}\label{classify}
Suppose $D_\Omega=\Omega+i\RR^4$ is a tube with affinely homogeneous base,
having as part of its boundary an affinely homogeneous non-degenerate
hypersurface. Then, up to biholomorphism, $D_\Omega$ must be exactly one of
the following:--
$$\begin{array}{rcl}
B_+^>&:=&\left\{z\in\CC^4:x_4>x_1^2+x_2^2+x_3^2\right\},\\
B_+^<&:=&\left\{z\in\CC^4:x_4<x_1^2+x_2^2+x_3^2\right\},\\
B_-^>&:=&\left\{z\in\CC^4:x_4>x_1^2+x_2^2-x_3^2\right\},\\
B_-^<&:=&\left\{z\in\CC^4:x_4<x_1^2+x_2^2-x_3^2\right\},\\
H^>  &:=&\left\{z\in\CC^4:x_4>x_1x_2+x_3^2\mbox{\rm\ and }
                       x_1>0\right\},\\
H^<  &:=&\left\{z\in\CC^4:x_4<x_1x_2+x_3^2\mbox{\rm\ and }
                        x_1>0\right\},\\
N_+^>&:=&\left\{z\in\CC^4:x_4>x_1x_2+x_3^2+x_1^2x_3+x_1^4\right\},\\
N_+^<&:=&\left\{z\in\CC^4:x_4<x_1x_2+x_3^2+x_1^2x_3+x_1^4\right\},\\
N_-^>&:=&\left\{z\in\CC^4:x_4>x_1x_2+x_3^2+x_1^2x_3-x_1^4\right\},\\
N_-^<&:=&\left\{z\in\CC^4:x_4<x_1x_2+x_3^2+x_1^2x_3-x_1^4\right\},\\
C^>  &:=&\left\{z\in\CC^4:x_4>x_1x_2+x_1x_3^2\mbox{\rm\ and }
                        x_1>0\right\},\\
C^<  &:=&\left\{z\in\CC^4:x_4<x_1x_2+x_1x_3^2\mbox{\rm\ and }
                        x_1>0\right\},\\
D^>  &:=&\left\{z\in\CC^4:x_4^2>x_1x_2+x_1^2x_3\mbox{\rm\ and }
                        x_1>0\right\},\\
D^<  &:=&\left\{z\in\CC^4:x_4^2<x_1x_2+x_1^2x_3\mbox{\rm\ and }
                              x_1>0\right\},
\end{array}$$
where $x_j=\hbox{Re}\,z_j$ for $j=1,\ldots,4$.
\end{theorem}
The domain $B_+^>$ is the ordinary unit ball written in tube form. The domains
$B_+^<$, $B_-^>$, $B_-^<$ are the tube realisations of domains lying to either
side of a quadric (\ref{quadric}) as mentioned earlier. They are the
(pseudo-) balls and, in the indefinite case, there are also affinely
homogeneous domains
$H^>$ and~$H^<$, which can be thought of as `half-pseudo-balls'. The
domains $N_+^>$, \dots, $N_-^<$ are nil-balls
in the sense of Penney~\cite{P1, P3}. This leaves $C^>$, $C^<$, $D^>$, $D^<$,
which are new. Remarkably, for each of these domains, the Levi non-degenerate
part of the boundary admits an explicit rational transformation to Chern-Moser
normal form. In this normal form the CR isotropy becomes linear. This enables
us to determine the full holomorphic automorphism group of the domains in
question.

The only pseudo-convex domain in Theorem~\ref{classify} is the unit
ball~$B_+^>$. Its complement, $B_+^<$ is pseudo-concave. The remaining domains
naturally fall into two classes
$$B_-^>,\,H^>,\,N_+^>,\,N_-^>,\,C^>,\,D^>\quad\mbox{versus}\quad
B_-^<,\,H^<,\,N_+^<,\,N_-^<,\,C^<,\,D^<$$
according to whether the Levi form of the non-degenerate part of the
boundary has type $++-$ or $+--$ (in some cases there is also a Levi
flat part to the
boundary).

The paper is organised as follows. In Section~\ref{exa} we examine the list of
\cite{EE2} and use it to prove Theorem~\ref{classify}. The new domains $D^>$
and $D^<$ are the main subject of the remainder of the article. We verify that
they do not have bounded realisations. In Section~\ref{autogr} we determine the
automorphism groups of these domains by means of exploiting the properties
of~$\Gamma$. These automorphism groups turn out to be Lie groups of
dimension~10, and their structure shows that the domains that $D^>$ and $D^<$
are not equivalent to any nil-balls. The new domains $C^>$ and $C^<$ are
similar but slightly easier. In Section~\ref{autoC} we sketch their analysis
along the lines already done in Section~\ref{autogr} for $D^>$ and~$D^<$.

We would like to acknowledge that this work
started while the third author was visiting the University of
Adelaide and continued during his visit to CIAM at the University
of South Australia. We would also like to thank Gerd Schmalz for
several useful conversations.

\section{Classification and New Examples}\label{exa}

The non-degenerate affinely homogeneous hypersurfaces with isotropy in $\RR^4$
are classified in~\cite{EE2}, where they are grouped according to their
complexification into 20 different types, some of which have more than one
real form. Sometimes, there is a parameter, as in hypersurface~\#4:--
$$
\{x\in\RR^4:x_4=x_1x_2+x_3^2+x_1^2x_3+\alpha x_1^4\},\quad\mbox{for }
\alpha\in\RR,
$$
where $x:=(x_1,x_2,x_3,x_4)$. There is a list of
explicit defining equations in \cite{EE2} and it is an
elementary matter to verify that each hypersurface $\Gamma$ from the list is,
indeed, affinely  homogeneous with isotropy. This is best done by infinitesimal
means: attached to any hypersurface is its symmetry algebra, namely the affine
vector fields on $\RR^4$ that are tangent to $\Gamma$ along~$\Gamma$, and the
question is whether this algebra is acting transitively near $\Gamma$'s chosen
basepoint and has dimension at least~$4$. This is mere computation and in many
cases, the symmetry algebra is exactly $4$-dimensional. In these cases, we may
manufacture a $4\times 4$ matrix from any chosen basis of affine vector fields
and compute its determinant: it will be a polynomial on~$\RR^4$. Where this
polynomial does not vanish, the symmetry algebra contains linearly independent
vector fields and so there are open orbits. Otherwise there cannot be open
orbits. All of these computations are easily carried out with computer algebra
and a Maple program, with the explicit defining equations already included, is
available\footnote{
ftp://ftp.maths.adelaide.edu.au/pure/meastwood/maple7/orbits}.
It turns out that this determinant is, in many cases, identically zero.
Hypersurfaces \#7--\#20 are thereby eliminated and we are left with (after some
minor affine changes of variables from~\cite{EE2}):--
\renewcommand{\arraystretch}{1.2}
\begin{center}\begin{tabular}{|l|l|c|}\cline{2-3}
\multicolumn{1}{c|}{}&Equation & Basepoint\\
\hline
\#1 &$x_4=x_1^2+x_2^2\pm x_3^2$               & $(0,0,0,0)$\\ \hline
\#2 &$x_1^2+x_2^2+x_3^3\pm x_4^2=1\quad\mbox{or}\quad
       x_1^2-x_2^2-x_3^3\pm x_4^2=1$           & $(1,0,0,0)$\\ \hline
\#3 &$x_4=x_1x_2+x_3^2+x_1^3$                 & $(0,0,0,0)$\\ \hline
\#4 &$x_4=x_1x_2+x_3^2+x_1^2x_3+\alpha x_1^4$ & $(0,0,0,0)$\\ \hline
\#5 &$x_4=x_1x_2+x_1x_3^2$                    & $(1,0,0,0)$\\ \hline
\#6 &$x_4^2=x_1x_2+x_1^2x_3$                  & $(1,0,1,1)$\\ \hline
\end{tabular}\end{center}\renewcommand{\arraystretch}{1}
\smallskip \noindent {\em Proof of Theorem~\ref{classify}:}\quad
We have just identified the only possible affinely homogeneous hypersurfaces
satisfying the hypotheses of Theorem~\ref{classify}. That is not to say that
all these possibilities actually occur and, in fact, \#2 does not. The affine
symmetry algebra of the sphere, for example, is spanned by the six vector
fields $x_j\partial/\partial x_k-x_k\partial/\partial x_j$ for $j<k$ and these
fields are nowhere linearly independent. In fact, the associated group is
${\mathrm{SO}}(4)$ with orbits concentric spheres and the origin. A similar
conclusion holds for other signatures.

All the remaining cases do, in fact, give rise to affinely homogeneous domains
as follows. Case~\#1, with either sign, has a $7$-dimensional symmetry algebra
${\mathcal{S}}$ and outside of $\Gamma$ itself these fields have maximal rank:
we obtain $B_+^>$, $B_+^<$, $B_-^>$, $B_-^<$. There is also the possibility of
having a Lie subalgebra ${\mathcal{R}}\subset{\mathcal{S}}$ with
$\dim{\mathcal{R}}=4,5,6$ and such that the fields in ${\mathcal{R}}$ have
maximal rank somewhere but not everywhere outside~$\Gamma$. This possibility is
easily investigated with computer algebra. First one finds the possible
subalgebras by parameterising the Grassmannians
${\mathrm{Gr}}_k({\mathcal{S}})$ for $k=4,5,6$ with affine coordinate patches
in the usual way, in each case solving the equations that these subspaces be
closed under Lie bracket. Then one computes the determinants of the $4\times 4$
minors for each of these subalgebras. Up to affine change of coordinates we
obtain just two possible affinely homogeneous `half-pseudo-balls'
$H^>$ and~$H^<$,
these domains having $5$-dimensional affine symmetry.

Case~\#3 has a $5$-dimensional symmetry algebra and the domains to either side
of $\Gamma$ are easily verified to be affinely homogeneous. It is a simple
computation to check that all of its $4$-dimensional subalgebras consist of
fields that are nowhere linearly independent. So, case~\#3 gives precisely the
following two holomorphically homogeneous domains
$$\{x\in\RR^4: x_4>x_1x_2+x_3^2+x_1^3\}+i\RR^4,\,
   \{x\in\RR^4: x_4<x_1x_2+x_3^2+x_1^3\}+i\RR^4.$$
However, the polynomial transformation
$$z_2 \mapsto z_2-3z_1^2/2\qquad z_4 \mapsto z_4-z_1^3/2$$
takes these two domains biholomorphically onto domains that we have already
listed, namely $B_-^>$ and $B_-^<$, respectively.

Case \#4 is the subject of \cite{EI} and will appear elsewhere. In brief, there
are three subcases according to whether $\alpha>1/12$, $\alpha=1/12$, or
$\alpha<1/12$. The case $\alpha=1/12$ again yields $B_-^>$ and $B_-^<$ by an
explicit polynomial change of coordinates. The equations of the other
hypersurfaces may be
placed in Chern-Moser normal form \cite{CM}:--
$$2\hbox{Im}\,w_4=w_1\overline{w_2}+w_2\overline{w_1}+|w_3|^2\pm|w_1|^4,$$
where the sign before $|w_1|^4$ is the sign of $\alpha-1/12$. This sign is a CR
invariant (as observed, for example, in~\cite{BSW}) and we obtain the domains
$N_+^>$, $N_+^>$, $N_-^>$, $N_-^<$. These were also found by
Penney \cite{P3} and are examples of his {\em nil-balls}.

In cases~\#4, \#5, and~\#6, the symmetry algebra is $4$-dimensional and so
there is no possibility of a Lie subalgebra giving rise to another domain. The
Maple computer program `orbits' mentioned early in this section finds the only
open orbits and these account for the remaining domains $C^>$, $C^<$, $D^>$,
$D^<$ in the statement of Theorem~\ref{classify}. That they are holomorphically
distinct and not nil-balls is proved in Section~\ref{autogr}, concerned with
$D^>$ and $D^<$, and Section~\ref{autoC}, concerned with $C^>$ and~$C^<$.
With these matters postponed, the theorem is proved. \qed

We conclude this section with some preliminary discussion of the cases~$D^>$
and $D^<$ before a detailed analysis of their holomorphic automorphisms in
Section~\ref{autogr}. So, consider the following affinely homogeneous
hypersurface in~$\RR^4$:--
$$\Gamma:=\left\{x\in\RR^4: x_4^2=x_1x_2+x_1^2x_3, x_1>0\right\},$$
and let
$$\begin{array}{rcl}
\Omega^{>}&:=&\left\{x\in\RR^4:x_4^2>x_1x_2+x_1^2x_3, x_1>0\right\},\\
\Omega^{<}&:=&\left\{x\in\RR^4:x_4^2<x_1x_2+x_1^2x_3, x_1>0\right\}.
\end{array}$$
Clearly, $\Gamma\subset\partial\Omega^{>}\cap \partial\Omega^{<}$. In fact,
these domains are not as similar as they might first appear. Part of the
boundary of $\Omega^>$ is Levi flat, namely
$\{x\in\RR^4: x_1=0\mbox{ and }x_4\not=0\}$. In the case of $\Omega^<$,
however, the inequality $x_1>0$
serves only to specify one of the two connected components of
$\{x\in\RR^4:x_4^2<x_1x_2+x_1^2x_3\}$.

One can verify that $\Omega^{>}$ and $\Omega^{<}$ are invariant under the group
$G$ that consists of the following affine transformations of $\RR^4$
$$
\begin{array}{lll}
x_1&\mapsto& q x_1,\\
x_2 &\mapsto& qr^2(s+t^2) x_1+qr^2 x_2+2qr^2t x_4,\\
x_3 &\mapsto& r^2 x_3-r^2s,\\
x_4 &\mapsto&qrt x_1+qr x_4,
\end{array}
$$
where $q>0$, $r\in\RR^*$, $s,t\in\RR$. We will now show that $G$ acts
transitively on each of $\Omega^{>}$ and $\Omega^{<}$. Take the point
$(1,0,0,1)\in\Omega^{>}$ and apply a mapping from $G$ to it. The result is the
point
$$
\Bigl(q, qr^2(t^2+2t+s),-r^2s,qr(t+1)\Bigr).
$$
Let $(x_1^0,x_2^0,x_3^0,x_4^0)$ be any other point in $\Omega^{>}$.
Then setting
$$
\begin{array}{cc}
\displaystyle q=x_1^0,&\displaystyle
r=\frac{\sqrt{(x_4^0)^2-x_1^0x_2^0-(x_1^0)^2x_3^0}}{x_1^0},\\ \\
\displaystyle t=\frac{x_4^0}{\sqrt{(x_4^0)^2-x_1^0x_2^0-(x_1^0)^2x_3^0}}-1,&
\displaystyle s=-\frac{(x_1^0)^2x_3^0}{(x_4^0)^2-x_1^0x_2^0-(x_1^0)^2x_3^0},
\end{array}
$$
we obtain an element of $G$ that maps $(1,0,0,1)$ into
$(x_1^0,x_2^0,x_3^0,x_4^0)$. This proves that $\Omega^{>}$ is affinely
homogeneous. A similar argument shows that $\Omega^{<}$ is affinely homogeneous
as well. It can be shown that $G$ is in fact the full group of affine
automorphisms of each of $\Omega^{>}$ and $\Omega^{<}$.

It then follows that the corresponding tube domains $D^{>}:=D_{\Omega^{>}}$ and
$D^{<}:=D_{\Omega^{<}}$ are holomorphically homogeneous since the group
$\tilde G$ generated by $G$ (viewed as a group of affine transformations of
$\CC^4$) and translations in the imaginary directions in $\CC^4$, acts
transitively on each of them.

We note that neither of these domains has a bounded realisation. In fact,
neither of them is Kobayashi-hyperbolic (we remark here that it is shown in
\cite{N} that any connected homogeneous Kobayashi-hyperbolic manifold is
biholomorphically equivalent to a bounded domain in complex space). Indeed,
the domain $D^{>}$ contains the affine complex line
$\{z\in\CC^4:z_1=1, z_2+z_3=0, z_4=1\}$ while $D^{<}$ contains the affine
complex line
$\{z\in\CC^4:z_1=1, z_2+z_3=1,z_4=0\}$.

\section{The Automorphism Groups of $D^{>}$ and $D^{<}$}\label{autogr}

Denote by $G^{>}$ and $G^{<}$ the groups of holomorphic automorphisms of
$D^{>}$ and $D^{<}$ respectively, equipped with the compact-open topology. In
this section we will determine these groups.

Let $\tilde\Gamma:=\Gamma+i\RR^4$. Clearly, $\tilde\Gamma$ is a connected Levi
non-degenerate hypersurface contained in $\partial D^{>}\cap\partial D^{<}$.
Since the Levi form of $\tilde\Gamma$ at every point has a positive and a
negative eigenvalue, every element of $G^{>}$ and $G^{<}$ extends past
$\tilde\Gamma$ to a biholomorphic mapping preserving $\tilde\Gamma$. Let $H$ be
the subgroup of $\tilde G$ given by the condition $r=1$. It is straightforward
to verify that $H$ acts simply transitively on $\tilde\Gamma$. Therefore, every
element $g$ of either $G^{>}$ or $G^{<}$, for a fixed $p\in\tilde\Gamma$, can
be uniquely represented as $g=h\circ t$, where $h\in H$ and $t$ is an element
of either $I_p^{>}$ or $I_p^{<}$, and $I_p^{>}$ and $I_p^{<}$ denote the
isotropy subgroups of $p$ in $G^{>}$ and $G^{<}$ respectively.

Thus, in order to determine $G^{>}$ and $G^{<}$ we must find $I_p^{>}$ and
$I_p^{<}$ for some $p\in\tilde\Gamma$. Let $p_0:=(1,0,1,1)$ and $I_{p_0}$ be
the group of all {\em local}\, holomorphic automorphisms of $\tilde\Gamma$
defined near $p_0$ and preserving $p_0$. Clearly, $I_{p_0}^{>}\subset I_{p_0}$
and $I_{p_0}^{<}\subset I_{p_0}$. Equipped with the topology of uniform
convergence of the derivatives of all orders of the component functions on
compact subsets of $p_0$, the group $I_{p_0}$ is known to carry the structure
of a real algebraic group (see, e.g., \cite{B}, \cite{L1}, \cite{VK}).

It is straightforward to show that $I_{p_0}\cap \tilde G$ consists of the
mappings
$$
\begin{array}{lll}
z_1 &\mapsto & z_1,\\
z_2 &\mapsto & 2(r^2-r)z_1+r^2z_2+2(r-r^2)z_4,\\
z_3 &\mapsto & r^2z_3+1-r^2,\\
z_4 &\mapsto & (1-r)z_1+rz_4,
\end{array}
$$
where $r\in\RR^*$. As the following proposition shows, $I_{p_0}$ is in fact
$3$-dimensional and admits an explicit description.

\begin{proposition}\label{descrisotrop} The group $I_{p_0}$ consists of all
mappings of the form
\begin{equation}
\begin{array}{lll}
z_1&\mapsto & z_1,\\
z_2&\mapsto & -2v^2z_1^3+i(-2v+4vr+u)z_1^2-4ivrz_1z_4\\
&&\quad+2(v^2+r^2-r)z_1+r^2z_2+2(r-r^2)z_4+2iv-iu,\\
z_3&\mapsto & 2v^2z_1^2-i(2u+4vr)z_1+r^2z_3+4ivrz_4-r^2\\
&&\quad+1-2v^2+2iu,\\
z_4&\mapsto &-ivz_1^2+(1-r)z_1+rz_4+iv,
\end{array}\label{formisotrop}
\end{equation}
where $r\in\RR^*$, $u,v\in\RR$.
\end{proposition}

\noindent {\em Proof:}\quad First, we will write the equation of $\tilde\Gamma$
near $p_0$ in Chern-Moser normal form \cite{CM}. Consider the change of
coordinates
$$\begin{array}{lll}
z_1 &=&\displaystyle \frac {11 w_1+4}{w_1+4},\\ \\
z_2 &=& \displaystyle -\frac {96i\left (4iw_2-5iw_3+11w_4
\right )}{5w_1+20}\\ \\
&&\displaystyle{}+\frac{1600i\left (4iw_2-5iw_3-6i
w_3^2+6w_4\right)}{\left(5w_1+20\right)^2}-
\frac{1280w_3^2}{\left(w_1+4\right)^3}+24 iw_4,\\ \\
z_3 &=& \displaystyle\frac{32i\left
(2\,iw_2+3w_4\right )}{5w_1+20}-\frac{32w_3^2}
{\left(w_1+4\right )^2}-4iw_4+1,\\ \\
z_4 &=& \displaystyle -\frac {8(6w_3+5)}{w_1+4}+
\frac{160w_3}{\left(w_1+4\right)^2}+11.
\end{array}$$
In the $w$-coordinates the point $p_0$ becomes the origin, and the equation of
$\tilde\Gamma$ near the origin takes the form
\begin{equation}
\hbox{Im}\,w_4=F\left(w',\overline{w'}\right):=
\frac{N\left(w',\overline{w'}\right)}{D\left(w',\overline{w'}\right)},
\label{neweq}
\end{equation}
where $w':=(w_1,w_2,w_3)$ and
$$\begin{array}{lll}
D\left(w',\overline{w'}\right)&=&
\left (11|w_1|^2+24 w_1+24\overline{w_1}+16\right)\left(5|w_1|^2+16\right),\\
N\left(w',\overline{w'}\right)&=&
8\hbox{Re}\,\Bigl(48w_1|w_3|^2+25w_1^2\overline{w_3}^2
+16|w_3|^2+36|w_1|^2|w_3|^2\\
&&\quad{}+2(11|w_1|^2+24w_1+24\overline{w_1}+16)w_1\overline{w_2}\Bigr).
\end{array}$$
Denote by $F_{k\overline{l}}$ the polynomial of degree $k$ in $w'$ and degree
$l$ in $\overline{w'}$ in the power series expansion of the function
$F\left(w',\overline{w'}\right)$ about the origin. A straightforward
calculation gives
\begin{equation}
\begin{array}l
F_{1\overline{1}}=
\displaystyle\frac{|w_3|^2}{2}+\hbox{Re}\,w_1\overline{w_2},\\ \\
F_{2\overline{2}}=
\displaystyle2\hbox{Re}\left(-\frac {5}{32}w_1^2\overline{w_1}\overline{w_2}
+\frac{25}{64}w_1^2\overline{w_3}^2+\frac{5}{16}|w_1|^2|w_3|^2\right),\\ \\
F_{3\overline{2}}=\displaystyle
-\frac {75}{128}w_1^3\overline{w_3}^2
-\frac {75}{64}w_1^2w_3\overline{w_1}\overline{w_3}-\frac
{75}{128}w_1w_3^2\overline{w_1}^2,\\ \\
F_{3\overline{3}}=\displaystyle2\hbox{Re}
\left(\frac{25}{512}w_1^3\overline{w_1}^2\overline{w_2}
+\frac{175}{128}w_1^3\overline{w_1}\overline{w_3}^2
+\frac{1425}{1024}|w_1|^4|w_3|^2\right).
\end{array}\label{impterms}
\end{equation}
It follows from the first formula in (\ref{impterms}) that the operator
$tr$ defined in \cite{CM} in this case is
$$tr= 2\left(\frac{\partial^2}{\partial w_1 \partial
\overline{w_2}} + \frac{\partial^2}{\partial w_2 \partial
\overline{w_1}} + \frac{\partial^2}{\partial w_3 \partial
\overline{w_3}} \right).$$
Equation~(\ref{neweq}) and the last three formulae in
(\ref{impterms}) now yield
$$F(w',0)=\frac{\partial F}{\partial \overline{w'}}(w',0)=
tr\,F_{2\overline{2}}=
tr^2\, F_{3\overline{2}}=tr^2\, F_{3\overline{3}}=0,$$
which shows that equation (\ref{neweq}) is indeed in Chern-Moser normal form.

We will now find $I_{p_0}$ in the $w$-coordinates. Let $I$ denote the group of
mappings of the form
$$\begin{array}{lll}
w_1&\mapsto& w_1,\\
w_2&\mapsto&
\displaystyle\left(i\mu-\frac{\nu^2}{2}\right)w_1+r^2w_2+i\nu rw_3,\\
w_3&\mapsto& i\nu w_1+rw_3,\\
w_4&\mapsto& r^2w_4,
\end{array}$$
where $r\in\RR^*$, $\mu,\nu\in\RR$. One can directly verify that $I$ is a
subgroup of $I_{p_0}$. We will show that in fact $I_{p_0}=I$.

Suppose that $I_{p_0}\ne I$ and let $f\in I_{p_0}\setminus I$. It then easily
follows from the condition that the term $F_{2\overline{2}}$ does not change
when $f$ is applied to equation (\ref{neweq}) that, without loss of generality,
the Jacobian matrix of $f$ at the origin can be assumed to have the form
\begin{equation}
\left (\begin{array}{cccc}
\beta&0&0&0\\
0&\beta&0&0\\
0&0&\beta&0\\
0&0&0&1
\end{array}\label{matrix}
\right),
\end{equation}
where $|\beta|=1$, $\beta\ne 1$. Since the local isotropy group of any Levi
non-degenerate hypersurface in $\CC^4$ is linearisable (see \cite{E1},
\cite{E2}), there exist normal coordinates in which transformations from $I$
remain linear and $f$ coincides with the linear transformation given by matrix
(\ref{matrix}). Since the term $F_{3\overline{2}}$ in these new coordinates has
to be preserved by $f$, it follows that $F_{3\overline{2}}=0$. On the other
hand, one can show that the term $F_{3\overline{2}}$ in formula
(\ref{impterms}) can not be eliminated by any transformation to the normal form
that leaves all mappings from $I$ linear. This follows from the fact that the
family of Chern-Moser chains invariant under the action of $I$ is characterised
by the condition that the tangent vector of a chain at the origin is a multiple
of a vector of the form $(0,a,0,1).$ Thus, for a normalisation $w\mapsto w^*$
that transforms such a chain into the line $\{w^{*'}=0,\,\hbox{Im}w_4^{*}=0\}$,
we have $\partial w^{*'}/\partial w_4(0) = (0,a,0,1)$. With this single
parameter freedom one can not eliminate the entire term $F_{3\overline{2}}$ in
formula (\ref{impterms}). This contradiction shows that $I_{p_0}=I$.

Recomputing the group $I_{p_0}$ in the $z$-coordinates, we obtain formulae
(\ref{formisotrop}) with $u=\frac{16}{25}\mu$ and $v=\frac{2}{5}\nu$. The
proposition is proved. \qed

Since all elements of $I_{p_0}$ are polynomial mappings, we have
$I_{p_0}\subset I_{p_0}^>$ and $I_{p_0}\subset I_{p_0}^<$, which implies that
$I_{p_0}^{>}=I_{p_0}^{<}=I_{p_0}$. Hence $G^{>}=G^{<}=H\circ I_{p_0}$. A
straightforward calculation now leads to the following

\begin{theorem}\label{descrautgroup} The automorphism groups of $D^{>}$ and
$D^{<}$ coincide and consist of all mappings of the form
\begin{equation}\label{descfullgroup}
\begin{array}{lll}
z_1&\mapsto& qz_1+i\alpha_1,\\
z_2&\mapsto& -2qv^2z_1^3+iq(-2v+4vr-2vt+u)z_1^2-4iqvrz_1z_4\\
&&\quad{}+q(2v^2+2r^2-2r-2tr+2t+s+t^2)z_1+qr^2z_2\\
&&\qquad{}+2q(r-r^2+tr)z_4+i\alpha_2,\\
z_3&\mapsto& 2v^2z_1^2-2i(2vr+u)z_1+r^2z_3+4ivrz_4\\
&&\quad{}-r^2-2v^2-s+1+i\alpha_3,\\
z_4&\mapsto& -iqvz_1^2+q(t-r+1)z_1+qrz_4+i\alpha_4,
\end{array}
\end{equation}
where $s,t,u,v,\alpha_1,\alpha_2,\alpha_3,\alpha_4\in\RR$, $q>0$, $r\in\RR^*$.
\end{theorem}

We will now show that neither of $D^{>}$, $D^{<}$ is equivalent to a nil-ball.
Otherwise one of these domains would be holomorphically equivalent to a domain
$D$ that has a Levi non-degenerate boundary and admits an action of a connected
nilpotent Lie group $N$ by holomorphic transformations defined in a
neighbourhood of $\overline{D}$, such that the induced action of $N$ on
$\partial D$ is transitive (see \cite{P1}, \cite{P3}). Let $F$ denote the
equivalence mapping. Since the Levi form of $\tilde\Gamma$ at every point has a
positive and a negative eigenvalue, $F$ extends to a neighbourhood of
$\tilde\Gamma$ to a biholomorphic mapping onto a neighbourhood of
$\overline{D}$, and we have $F(\tilde\Gamma)=\partial D$ (the values of $F$
near the part of the boundary intersecting the hyperplane
$\{\hbox{Re}\,z_1=0\}$
necessarily approach infinity). Therefore, $N$ acts on either $D^{>}$ or
$D^{<}$ by biholomorphic transformations defined in a neighbourhood of
$\tilde\Gamma$ in such a way that the induced action of $N$ on $\tilde\Gamma$
is transitive. This implies that the Lie algebra ${\mathfrak g}$ of
$G^{>}=G^{<}$ has a nilpotent subalgebra that acts transitively on
$\tilde\Gamma$. We will now show that such a subalgebra in fact does not exist.

The algebra ${\mathfrak g}$ is spanned by the following holomorphic vector
fields
$$\begin{array}{l}
Z_1:= \displaystyle z_1\frac{\partial}{\partial z_1}
+z_2\frac{\partial}{\partial z_2}+z_4\frac{\partial}{\partial z_4},\,\,
Z_2:= \displaystyle z_1\frac{\partial}{\partial z_2}
-\frac{\partial}{\partial z_3},\\ \\
Z_3:=\displaystyle 2z_4\frac{\partial}{\partial z_2}
+z_1\frac{\partial}{\partial z_4},\,\,
Z_4:=\displaystyle i\frac{\partial}{\partial z_1}
+i\frac{\partial}{\partial z_4},\\ \\
Z_5:=\displaystyle i\frac{\partial}{\partial z_2},\,\,
Z_6:=\displaystyle i\frac{\partial}{\partial z_3},\,\,
Z_7:=\displaystyle 2i\frac{\partial}{\partial z_2}
+i\frac{\partial}{\partial z_4},\\ \\
Z_8:=\displaystyle 2(z_1+z_2-z_4)\frac{\partial}{\partial z_2}
+2(z_3-1)
\frac{\partial}{\partial z_3}+(z_4-z_1)\frac{\partial}{\partial z_4},\\ \\
Z_9:=\displaystyle iz_1^2\frac{\partial}{\partial z_2}
-2iz_1\frac{\partial}{\partial z_3},\\ \\
Z_{10}:=\displaystyle 2i(z_1^2-2z_1z_4)\frac{\partial}{\partial z_2}
-4i(z_1-z_4)\frac{\partial}{\partial z_3}
-iz_1^2\frac{\partial}{\partial z_4}.
\end{array}$$
The commutation relations for the vector fields above are as follows
\begin{center}\begin{tabular}{c|c|c|c|c|c|c|c|c|c|c}
$[Z_i,Z_j]$
&$Z_1$&$Z_2$&$Z_3$&$Z_4$&$Z_5$&$Z_6$&$Z_7$&$Z_8$&$Z_9$&$Z_{10}$\\ \hline
$Z_1$&&$0$&$0$&$-Z_4$&$-Z_5$&$0$&$-Z_7$&$0$&$Z_9$&$Z_{10}$\\ \hline
$Z_2$&&&$0$&$-Z_5$&$0$&$0$&$0$&$2Z_2$&$0$&$0$\\ \hline
$Z_3$&&&&$-Z_7$&$0$&$0$&$-2Z_5$&$Z_3$&$0$&$-2Z_9$\\ \hline
$Z_4$&&&&&$0$&$0$&$0$&$0$&$-2Z_2$&$2Z_3$\\ \hline
$Z_5$&&&&&&$0$&$0$&$2Z_5$&$0$&$0$\\ \hline
$Z_6$&&&&&&&$0$&$2Z_6$&$0$&$0$\\ \hline
$Z_7$&&&&&&&&$Z_7$&$0$&$4Z_2$\\ \hline
$Z_8$&&&&&&&&&$-2Z_9$&$-Z_{10}$\\ \hline
$Z_9$&&&&&&&&&&$0$\\ \hline
$Z_{10}$&&&&&&&&&&
\end{tabular}\end{center}
If a subalgebra ${\mathfrak h}\subset{\mathfrak g}$ acts transitively on
$\tilde\Gamma$, it must contain vector fields $Z_1':=Z-Z_1$, $Z_4':=Z_4+W$,
with some $Z,W\in{\mathfrak i}$, where ${\mathfrak i}$ is the Lie algebra of
$I_{p_0}$ (observe that ${\mathfrak i}$ is spanned by $Z_8$, $Z_9-Z_5+2Z_6$,
$Z_{10}+Z_7$). It follows from the commutation relations above that
$[Z_1',Z_4']$ is a vector field of the form $Z_4+L$, where $L$ is a linear
combination of $Z_2$, $Z_3$, $Z_5$, $Z_6$, $Z_7$, $Z_8$, $Z_9$, $Z_{10}$.
Thus, by induction, all the higher-order commutators
$[Z_1',[Z_1',\dots,[Z_1', [Z_1',Z_4']]\dots]]$ have this form and hence are
non-zero. Therefore, ${\mathfrak h}$ is not nilpotent.

Thus, we have proved the following:--
\begin{theorem}
Each of the domains $D^{>}, D^{<}\subset\CC^4$ is holomorphically homogeneous,
does not have a bounded realisation, and is not equivalent to a nil-ball.
\end{theorem}

\section{The Automorphism Groups of $C^{>}$ and $C^{<}$}\label{autoC}
The discussion for $C^>$ and $C^<$ largely parallels that for $D^>$ and $D^<$
just given in the previous section. In this section we omit all details.

Is is easily verified that affine transformations of the form
\begin{equation}\label{afftrans}\begin{array}{rcl}
x_1&\mapsto& qx_1,\\
x_2&\mapsto& r^2(x_2-2sx_3+t),\\
x_3&\mapsto& r(x_3+s),\\
x_4&\mapsto& qr^2(x_4+s^2x_1+tx_1),\end{array}\end{equation}
where $q>0$, $r\in\RR^*$, $s,t\in\RR$,
form a group and that this group acts transitively on
$$\{x\in\RR^4:x_4>x_1x_2+x_1x_3^2\mbox{ and }x_1>0\}$$
and
$$\{x\in\RR^4:x_4<x_1x_2+x_1x_3^2\mbox{ and }x_1>0\}.$$
It follows that $C^>$ and $C^<$ are holomorphically homogeneous.
Neither has a bounded realisation since the complex lines
$\{z\in\CC^4:z_2=z_3=0,z_4=\pm1\}$ lie entirely in $C^>$ and $C^<$,
respectively.

As in Section~\ref{autogr}, to determine the full holomorphic automorphism
groups of these domains, it suffices to know the local CR isotropy of the real
Levi-indefinite hypersurface
$$\tilde\Gamma:=\{z\in\CC^4:x_4=x_1x_2+x_1x_3^2\mbox{ and }x_1>0\}$$
near the basepoint $(1,0,0,0)$. Already we have an isotropy transformation
coming from the affine transformations (\ref{afftrans}) of~$\RR^4$, namely
$$
\begin{array}{lll}
z_1&\mapsto& z_1,\\
z_2&\mapsto& r^2z_2,\\
z_3&\mapsto& rz_3,\\
z_4&\mapsto& r^2z_4,
\end{array}
$$
for $r\in\RR^*$. We claim that the CR isotropy is generated, in addition, by
$$
\begin{array}{lll}
z_1&\mapsto& z_1,\\
z_2&\mapsto& z_2+iu(z_1-1),\\
z_3&\mapsto& z_3,\\
z_4&\mapsto& \displaystyle z_4+\frac{iu(z_1^2-1)}{2}
\end{array}
$$
and
\begin{equation}
\begin{array}{lll}
z_1&\mapsto& z_1,\\
z_2&\mapsto& \displaystyle z_2+\frac{(e^{2i\theta}-1)z_3^2}{2},\\
z_3&\mapsto& e^{i\theta}z_3,\\
z_4&\mapsto& z_4,
\end{array}\label{circleaction}
\end{equation}
where $u,\theta\in\RR$. To see this, we check that the holomorphic
change of variables
$$
\begin{array}{lll}
z_1 &=& w_1+1,\\
z_2 &=& \displaystyle w_2-\frac{iw_1w_4}{10}-\frac{2w_3^2}{(w_1+2)^2},\\
z_3 &=& \displaystyle \frac{2w_3}{w_1+2},\\
z_4 &=& \displaystyle
-iw_4+w_2+\frac{w_1(10w_2-i(w_1+2)w_4)}{20}
\end{array}
$$
takes the equation of the surface into Chern-Moser normal form:--
$$\hbox{Im}\,w_4=
10\hbox{Re}\,\frac{
4w_3\overline{w_3}(1+w_1)+2w_2w_1\overline{w_1}
+\overline{w_2}w_1^2\overline{w_1}+4\overline{w_2}w_1+2\overline{w_2}w_1^2}
{(2+w_1)(2+\overline{w_1})(20-w_1\overline{w_1})}$$
and that, in these coordinates, the CR isotropy is linear:--
$$
\begin{array}{lll}
w_1&\mapsto& w_1,\\
w_2 &\mapsto& r^2(w_2+iuw_1),\\
w_3 &\mapsto& re^{i\theta}w_3,\\
w_4 &\mapsto& r^2w_4,
\end{array}
$$
for $r\in\RR^*$, $u,\theta\in\RR$. In the original coordinates we obtain
precisely the stated isotropy. With the remaining affine symmetries
((\ref{afftrans}) with $r=1$) and imaginary translations, we have now
found the full holomorphic automorphism group. As a convenient basis for the
corresponding Lie algebra of holomorphic vector fields we may take
$$\begin{array}{l}
Z_1:= \displaystyle
z_1\frac{\partial}{\partial z_1}+z_4\frac{\partial}{\partial z_4},\,\,
Z_2:= \displaystyle
-2z_3\frac{\partial}{\partial z_2}+\frac{\partial}{\partial z_3},\,\,
Z_3:=\displaystyle
\frac{\partial}{\partial z_2}+z_1\frac{\partial}{\partial z_4},\,\,\\ \\
Z_4:=\displaystyle
i\frac{\partial}{\partial z_1},\,\,
Z_5:=\displaystyle
i\frac{\partial}{\partial z_2},\,\,
Z_6:=\displaystyle
i\frac{\partial}{\partial z_3},\,\,
Z_7:=\displaystyle
i\frac{\partial}{\partial z_4},\\ \\
Z_8:=\displaystyle
2z_2\frac{\partial}{\partial z_2}
+z_3\frac{\partial}{\partial z_3}+2z_4\frac{\partial}{\partial z_4},\,\,
Z_9:=\displaystyle
2iz_1\frac{\partial}{\partial z_2}+iz_1^2\frac{\partial}{\partial z_4},\\ \\
Z_{10}:=\displaystyle
-iz_3^2\frac{\partial}{\partial z_2}
+iz_3\frac{\partial}{\partial z_3}
\end{array}$$
with commutation relations
\begin{center}\begin{tabular}{c|c|c|c|c|c|c|c|c|c|c}
$[Z_i,Z_j]$
&$Z_1$&$Z_2$&$Z_3$&$Z_4$&$Z_5$&$Z_6$&$Z_7$&$Z_8$&$Z_9$&$Z_{10}$\\ \hline
$Z_1$&&$0$&$0$&$-Z_4$&$0$&$0$&$-Z_7$&$0$&$Z_9$&$0$\\ \hline
$Z_2$&&&$0$&$0$&$0$&$2Z_5$&$0$&$Z_2$&$0$&$Z_6$\\ \hline
$Z_3$&&&&$-Z_7$&$0$&$0$&$0$&$2Z_3$&$0$&$0$\\ \hline
$Z_4$&&&&&$0$&$0$&$0$&$0$&$-2Z_3$&$0$\\ \hline
$Z_5$&&&&&&$0$&$0$&$2Z_5$&$0$&$0$\\ \hline
$Z_6$&&&&&&&$0$&$Z_6$&$0$&$-Z_2$\\ \hline
$Z_7$&&&&&&&&$2Z_7$&$0$&$0$\\ \hline
$Z_8$&&&&&&&&&$-2Z_9$&$0$\\ \hline
$Z_9$&&&&&&&&&&$0$\\ \hline
$Z_{10}$&&&&&&&&&&
\end{tabular}\end{center}

The rest of the discussion follows exactly the corresponding discussion in
Section~\ref{autogr} (except that ${\mathfrak{i}}$ is spanned by
$Z_8$, $Z_9-2Z_5-Z_7$, $Z_{10}$) and we have proved:--
\begin{theorem}
Each of the domains $C^{>}, C^{<}\subset\CC^4$ is holomorphically homogeneous,
does not have a bounded realisation, and is not equivalent to a nil-ball.
\end{theorem}
Finally notice that $C^>$ and $C^<$ are holomorphically distinguished from
$D^>$ and $D^<$, respectively, by their admitting a holomorphic circle
action~(\ref{circleaction}).

\end{document}